\documentclass[a4paper,11pt]{article}
\usepackage{graphicx}
\usepackage{amssymb}
\usepackage{latexsym, bm}
\usepackage{multicol}
\usepackage{indentfirst}
\usepackage{amssymb,amsfonts}
\usepackage{amsmath}
\usepackage{mathrsfs}
\newtheorem{theorem}{Theorem}
\newtheorem{lemma}{Lemma}
\newtheorem{corollary}{Corollary}

\newtheorem{definition}{Definition}

\newtheorem{fact}{Fact}
\usepackage{amssymb}
\title{Sparse Multipartite Graphs as Partition Universal  \\ for Graphs of Bounded-Degrees \thanks{Supported in part by NSFC and CSC.}}
\author{Qizhong Lin$^1$,\,\,Yusheng Li$^2$  \vspace*{0.3cm}\\
$^1$Center for Discrete Mathematics, Fuzhou University \\Fuzhou 350108, China\\
$^2$Department of Mathematics, Tongji University  \\Shanghai 200092, China \\
{\small\em Email:  linqizhong@fzu.edu.cn, li\_yusheng@tongji.edu.cn }}
\date{}
\begin{document}
\maketitle
\begin{abstract}
For graphs $G$ and $H$, let $G\to (H,H)$ signify that any red/blue edge coloring of $G$ contains a monochromatic $H$ as a subgraph,
and $\mathcal{H}(\Delta,n)=\{H:|V(H)|=n,\Delta(H)\le \Delta\}$.
For fixed $\Delta$ and $n$, we say that $G$ is a partition universal graph for $\mathcal{H}(\Delta,n)$ if $G\to (H,H)$ for every $H\in\mathcal{H}(\Delta,n)$.

In 1983, Chv\'{a}tal, R\"{o}dl, Szemer\'{e}di and Trotter proved that
for any $\Delta\ge2$ there exists a constant $B$ such that, for any $n$, if $N\ge Bn$ then $K_N$ is partition universal for $\mathcal{H}(\Delta,n)$.
Recently, Kohayakawa, R\"{o}dl, Schacht and Szemer\'{e}di proved that the complete graph $K_N$ in above result can be replaced by sparse graphs.
They obtained that for fixed $\Delta\ge2$, there exist constants $B$ and $C$ such that if $N\ge Bn$ and $p=C(\log N/N)^{1/\Delta}$,
then {\bf a.a.s.} $G(N,p)$ is partition universal graph for $\mathcal{H}(\Delta,n)$, where $G(N,p)$ is the standard  random graph on $N$ vertices with $\mathbb{P}(e)=p$ for each edge $e$. From some results of Bollob\'as and {\L}uczak, we know that {\bf a.a.s.} $\chi(G(N,p)) = \Theta((N/\log N)^{1-1/\Delta})$.

In this paper, we show that the random graph $G(N,p)$ in the above result can be replaced by the random multipartite graph.
Let $K_{r}(N)$ be the complete $r$-partite graph with $N$ vertices in each part, and $G_r(N,p)$ the random spanning subgraph of $K_r(N)$, in which each edge appears with probability $p$. It is shown that for fixed $\Delta\ge2$ there exist constants $r, B$ and $C$ depending only on $\Delta$ such that if $N\ge Bn$ and $p=C(\log N/N)^{1/\Delta}$, then {\bf a.a.s.} $G_r(N,p)$ is partition universal graph for $\mathcal{H}(\Delta,n)$.
The proof mainly uses the sparse multipartite regularity lemma.

\medskip

\end{abstract}

\section{Introduction}

For a family $\cal{H}$ of graphs, a graph $G$ is said to be an $\cal{H}$-{\em universal} if it contains $H$ as a subgraph for each $H\in\cal{H}$.
The construction of sparse universal graphs for various families
arises in the study of VLSI circuit design, and has attracted much of attention. VLSI (Very-large-scale integration) is the process of creating an integrated circuit by combining thousands of transistors into a single chip. We always view the circuits as graphs, and wish to deal with families of graphs. It is pointed out in \cite{chung-rose83} that the problem of designing an efficient single circuit, specialized for a variety of other circuits, can be viewed as constructing a small universal graph.
For the references of universal graphs, one can see, e.g.,
\cite{alon-asodi02,alon-cap07,alon-cap08,alon-cap-ko00,alon-cap-ko01,alon-kri-07,alstrup02,babai-chung-gr-sp82,Balogh10,bhatt-chung89,
capalbo99,capalbo-kosa99,chung90,chung-graham78,chung-graham79,chung-graham83,chung-graham-pi78,della08,friedman-pip87,kim2014,rodl81} and their references.

Let $G(n,p)$ be the standard  random graph on $n$ vertices with $\mathbb{P}(e)=p$ for each edge $e$. Let $K_{r}(n)$ be the complete $r$-partite graph with $n$ vertices in each part, and $G_r(n,p)$  the random spanning subgraph of $K_r(n)$ with $\mathbb{P}(e)=p$ for each edge $e$ of $K_r(n)$. Let us set
\[
\mathcal{H}(\Delta,n)=\{H\subseteq K_n: \Delta(H)\le \Delta\}, \hspace{3mm}\mbox{and}\hspace{3mm}
\mathcal{H}(\Delta,n,n)=\{H\subseteq K_{n,n}: \Delta(H)\le \Delta\},
\]
where $H\subseteq K_n$ and $H\subseteq K_{n,n}$ mean that $H$ is a spanning subgraph of $K_n$ and $K_{n,n}$, respectively. By counting all unlabeled $\Delta$-regular graphs on $n$ vertices, it is shown \cite{alon-cap-ko00} that each $\mathcal{H}(\Delta,n)$-universal graph has at least
$\Theta\big(n^{2-2/\Delta}\big)$ edges, where and henceforth $\Theta(f(n))$ to signify a function that differs from $f(n)$ up to a multiplicative positive constant. This lower bound was almost matched that from  \cite{alon-cap07, alon-cap08},
which constructed  $\mathcal{H}(\Delta,n)$-universal graphs of order $n$ with at most $\Theta(n^{2-2/\Delta}(\log n)^{4/\Delta})$ edges.
It is shown that {\bf a.a.s.} $G(n,p)$ is $\mathcal{H}(\Delta,n)$-universal with $p=c((\log n)/n)^{1/\Delta}$ in \cite{della-koha14},
and {\bf a.a.s.} $G_2(n,p)$ is  $\mathcal{H}(\Delta,n,n)$-universal with $p=c((\log n)/n)^{1/\Delta}$ in \cite{della-koha08,della-koha12}.
Here we say $G_{n,p}$ possesses a property ${\cal P}$ asymptotically almost surely ({\bf a.a.s.}) if $\mathbb{P}[G_{n,p}\in \mathcal{P}] = 1-o(1)$ as $n\to\infty$.

For graphs $G$ and $H$, let $G \to (H,H)$ signify that any red/blue edge coloring of $G$ contains a monochromatic $H$ as a subgraph.
Thus the Ramsey number $R(H)=\min\{N:\;K_N\to (H,H)\}$.
Furthermore, we say that a graph $G$ is {\em partition universal} for $\mathcal{H}(\Delta,n)$ if $G \to (H,H)$ for each $H\in\mathcal{H}(\Delta,n)$.

A well-known result of Chv\'atal, R\"odl, Szemer\'edi and Trotter \cite{ch-ro-sz-tr} is as follows,
which implies that $R(H)\le \Theta(n)$ for $H\in{\cal H}(\Delta,n)$.
The proof is a remarkable application of Szemer\'edi regularity lemma \cite{sze78}, in which they used the general form of the lemma.
\begin{theorem} [\cite{ch-ro-sz-tr}]\label{r-linear}
For any fixed integer $\Delta\ge2$, there exists a constant $B=B(\Delta)$ such that if $N\ge Bn$ then $K_N$ is partition universal for ${\cal H}(\Delta,n)$.
\end{theorem}

The size Ramsey number $\widehat{r}(H)$ is defined to be $\min\{e(G): G\to(H,H)\}$ in \cite{efrs-78-size}, where $e(G)$ is the size of the edge set of $G$.
Recently, by using the sparse regularity lemma (see \cite{koha97,koha03}), an elegant result of Kohayakawa, R\"{o}dl, Schacht and Szemer\'{e}di \cite{Kohayakawa} strengthens Theorem \ref{r-linear} by replacing $K_N$ in Theorem \ref{r-linear} with sparse sparse graphs, implying $\widehat{r}(H) \le \Theta(n^{2-1/\Delta}(\log n)^{1/\Delta})$ for $H\in{\cal H}(\Delta,n)$, and confirming the upper bound in a conjecture of R\"{o}dl and Szemer\'{e}di \cite{rodl-sze}.
\begin{theorem}[\cite{Kohayakawa}] \label{koh}
For fixed $\Delta\ge2$, there exist constants $B=B(\Delta)$ and $C=C(\Delta)$ such that if $N\ge Bn$ and $p=C(\log N/N)^{1/\Delta}$, then
\[
\lim_{n\to\infty}\mathbb{P}\big(G(N,p)\,\,\mbox{is partition universal for}\,\, {\cal H}(\Delta,n)\big)=1.
\]
\end{theorem}

From some well-known results of  Bollob\'as \cite{bollo88} and {\L}uczak \cite{luczak91} on chromatic numbers of random graphs with $p=C(\log N/N)^{1/\Delta}$ as above, we have $\chi(G(N,p)) = \Theta((N/\log N)^{1-1/\Delta})$.
In this paper, we shall show a version of Theorem \ref{koh} on sparse multipartite graphs as follows.

\begin{theorem}\label{main}
For fixed $\Delta\ge2$, there exist constants $r=r(\Delta), B=B(\Delta)$ and $C=C(\Delta)$ such that if $N\ge Bn$ and $p=C(\log N/N)^{1/\Delta}$, then
\[
\lim_{n\to\infty}\mathbb{P}\big(G_r(N,p)\,\,\mbox{is partition universal for}\,\, {\cal H}(\Delta,n)\big)=1.
\]
\end{theorem}

{\em Remark.} Note that the constant $r=r(\Delta)$ in Theorem \ref{main} depends on $\Delta$ only.
The main strategy for the proof of Theorem \ref{main} is the sparse multipartite regularity lemma \cite{koha97, koha03}. A key tool for the proof is an elegant embedding method developed in \cite{Kohayakawa}.
Moreover, the above result also holds if we color the edges of $G_r(N,p)$ by $k\ge 3$ colors, for which the constants $r$, $B$ and $C$ will depend on both $\Delta$ and $k$.

\section{Sparse multipartite regularity lemma}

In this section, we shall restate a sparse $r$-partite regularity lemma and some related properties.
In the next of this paper, let $K_{r}(N)$ be the complete $r$-partite  graph on parts $V^{(1)},V^{(2)},\dots,V^{(r)}$ with $|V^{(i)}|=N$ for $1\le i\le r$, and let $G_r=G_r(N)=G_r(V^{(1)},V^{(2)},\dots,V^{(r)})$ be the spanning subgraph of $K_{r}(N)$.
Suppose $0 < p,\,\eta \le 1$. For subsets $X\subseteq V^{(i)}$ and $Y\subseteq V^{(j)}$, $i\neq j$,
let $e_{G_r}(X,Y)$ be the number of edges of $G_r$ between $X$ and $Y$, and let
\[
d_{G_r,p}(X,Y)=\frac{e_{G_r}(X,Y)}{p|X||Y|},
\]
which is referred to as the $p$-{\em{density}} of the pair $(X,Y)$. We say that $G_r$ is an $(\eta,\lambda)$-bounded graph with respect to $p$-density if  any pair of sets $X\subseteq V^{(i)}$ and $Y\subseteq V^{(j)}$ $(i\neq j)$ with $|X|\ge \eta|V^{(i)}|,|Y|\ge\eta|V^{(j)}|$ satisfy
\[
d_{G_r,p}(X,Y) \le \lambda.
\]
For fixed $\epsilon>0$, we say such a pair $(X,Y)$ is $(\epsilon,p)$-{\em regular} if for all $X'\subseteq X$ and $Y'\subseteq Y$ with
\[
|X'|\ge \epsilon |X| \,\,\mbox{and}\,\, |Y'|\ge \epsilon |Y|,
\]
we have
\[
|d_{G_r,p}(X,Y)-d_{G_r,p}(X',Y')|\le\epsilon.
\]
Note that for $p = 1$ we get the well-known definition of $\epsilon$-regularity. The following is a variant
of the Szemer\'{e}di's regularity lemma \cite{sze78} for sparse multipartite graphs, developed independently by Kohayakawa and R\"{o}dl, see \cite{koha97,koha03}.
We restated it as follows.
\begin{theorem}[Sparse multipartite regularity lemma]\label{sparse-reg}
For any fixed $\epsilon>0$, $\lambda>1$, $t_0\ge 1$ and $r\ge 2$, there exist $T_0$, $\eta$ and $N_0$,
such that each $r$-partite graph $G_r(N)=G_r(V^{(1)},V^{(2)},\dots,V^{(r)})$ with $N\ge N_0$  that is $(\eta,\lambda)$-bounded with respect to density $p$ with  $0<p\le1$, has a partition
$\{V^{(i)}_0,V^{(i)}_1,\cdots, V^{(i)}_t\}$ for each $V^{(i)}$, where $t$ is same for each part $V^{(i)}$ and $t_0\le t \le T_0$, such that
\begin{align*}
(i)\,\, &|V^{(i)}_0|\le\epsilon N\,\,\mbox{and}\,\,|V^{(i)}_1|=|V^{(i)}_2|=\dots=|V^{(i)}_t|\,\,\mbox{for}\,\, 1\le i\le r;
\\(ii)\,\, & \mbox{All but at most}\,\,\epsilon t^2\,\,\mbox{pairs}\,\,\big(V^{(i)}_a,V^{(j)}_b\big)\,\,\mbox{are}\,\, (\epsilon,p)-regular\,\,\mbox{for}\,\,1\le i\neq j\le r\,\,\mbox{and}\,\, 1\le a,b \le t.
\end{align*}
\end{theorem}

In the following, let us introduce the hereditary nature of sparse regularity which specially holds for $r$-partite graphs, see \cite{gerke07}.
\begin{definition}
Let $\alpha,\epsilon>0$, and $0<p \le1$ be given and let $G_r(N)$ be defined above. For sets $X\subseteq V^{(i)}$ and $Y\subseteq V^{(j)}$ ($i\neq j$), we say that the pair
$(X,Y)$ is $(\epsilon,\alpha,p)$-dense if for all $X'\subseteq X$ and $Y'\subseteq Y$ with $|X'|\ge\epsilon|X|$ and $|Y'|\ge\epsilon|Y|$, we have
\[
d_{G_r,p}(X',Y')\ge\alpha-\epsilon.
\]
\end{definition}
It follows immediately from the definition that $(\epsilon,\alpha,p)$-denseness is inherited on large sets, i.e., that for an $(\epsilon,\alpha,p)$-dense pair
$(X,Y)$ with $|X'|\ge\mu|X|$ and $|Y'|\ge\mu|Y|$ the pair $(X',Y')$ is $(\epsilon/\mu,\alpha,p)$-dense. The
following result from \cite{gerke07} states that this ``denseness-property" is even inherited on
randomly chosen subsets of much smaller size with overwhelming probability.


\begin{theorem}[\cite{gerke07}, Corollary 3.8] \label{num-den-pair}
For every $\alpha,\beta>0$ and $\epsilon'>0$, there exist $\epsilon_0=\epsilon_0(\alpha,\beta,\epsilon')>0$ and $L=L(\alpha,\epsilon')$ such that,
 for any $0<\epsilon\le\epsilon_0$ and $0<p<1$, every $(\epsilon,\alpha,p)$-dense pair $(X,Y)$  in a graph $G$ has the following property:
 the number of pairs $(X',Y')$ with $X'\subseteq X$ and $Y'\subseteq Y$ with $|X'|=w_1\ge L/p$ and $|Y'|=w_2\ge L/p$ such that
$(X',Y')$ is an $(\epsilon',\alpha,p)$-dense pair is at least $(1-\beta^{\min\{w_1,w_2\}}){|X|\choose w_1}{|Y|\choose w_2}$.
\end{theorem}

\section{Properties of the random multipartite graph}

In this section, we dedicate in establishing the properties for the random multipartite graph $G_r(N,p)$.
Recall that $K_{r}(N)$ has $r$ parts $V^{(1)},V^{(2)},\dots,V^{(r)}$ with $|V^{(i)}|=N$ for $1\le i\le r$,
and $G_r(N,p)$ is the random $r$-partite spanning subgraph of $K_r(N)$ with probability $p$ for edge  appearance.


\begin{definition} \label{e(UiUj)}
For an integer $N$ and $0<p\le 1$, we say that a graph $G_r(N)$ has the property ${\mathscr{U}}_{N,p}$ if
for $1\le i\neq j\le r$, all $U^{(i)}\subseteq V^{(i)},U^{(j)}\subseteq V^{(j)}$
with $|U^{(i)}|,|U^{(j)}|\ge {N}/{\log N}$ satisfy
\[
|d_{G_r(N,p)}(U^{(i)},U^{(j)})-1|\le\frac{1}{\log N}.
\]
\end{definition}

The Chernoff's inequality (see e.g., \cite{cher,alo-spe,janson}) will be useful for the proof of the following fact.
\begin{lemma}\label{chernoff}
Let $X$ be a binomial random variable.
 If $0<\delta\le1$, then
 \[
 \mathbb{P}\big(|X-\mathbb{E}(X)|\ge\delta\mathbb{E}(X)\big)\le2\exp\big(-{\delta^2}\mathbb{E}(X)/3\big).
 \]
\end{lemma}

The first fact we will obtain as follows implies that  {\bf a.a.s.} the edge number for all sufficiently large pair defined as above are concentrated on the expectation for suitable probability $p$.
\begin{fact}\label{uniform}
If $p\ge12(\log N)^4/N$, then $\mathbb{P}(G_r(N,p)\in {\mathscr{U}}_{N,p})\to 1$  as $N\to\infty$.
\end{fact}

\noindent{\bf Proof.} The graph $G_r(N,p)\not\in{\mathscr{U}}_{N,p}$ means that there exists a pair $(U^{(i)},U^{(j)})$ (for some $1\le i\neq j\le r$) such that $|d_{G_r(N,p)}(U^{(i)},U^{(j)})-1|>\frac{1}{\log N}$, i.e.,
$$\left|e_{G_r(N,p)}(U^{(i)},U^{(j)})- p|U^{(i)}||U^{(j)}|\right|> p|U^{(i)}||U^{(j)}|/\log N.$$
 Hence, Lemma \ref{chernoff} implies the probability that the graph $G_r(N,p)\not\in{\mathscr{U}}_{N,p}$ is at most
\begin{align*}
{r\choose2}&{N\choose |U^{(i)}|}{N\choose |U^{(j)}|}\times2\exp\left(\frac{-p|U^{(i)}||U^{(j)}|}{3\log^2 N}\right)
\\&\le {r^2}\cdot\exp\left(|U^{(i)}|\log N+|U^{(j)}|\log N-\frac{p|U^{(i)}||U^{(j)}|}{3\log^2 N}\right),
\end{align*}
which will tends to zero as $N\to\infty$ since $p\ge12(\log N)^4/N$. \hfill$\Box$

\smallskip
For a $r$-partite graph $G_r(N)$ and integers $1\le k\le\ell< r$, let $\mathcal{K}$ be the family consists of all $k$-subset $K$ of $\bigcup_{i=1}^\ell V^{(i)}$ such that $|K\cap V^{(i)}|\le1$ for $1\le i\le \ell$.  Define the auxiliary bipartite graph $\Gamma=\Gamma(k,G_r(N))$ with color classes $\mathcal{K}$ and $\bigcup_{i=1}^\ell V^{(i)}$, where $(K,v)\in E(\Gamma)$  if and only if $\{w,v\}\in E(G_r(N))$ for all $w\in K$. Here $E(\Gamma)$ is the edge set of $\Gamma$.

\begin{definition}\label{congestion-def}
 Let integers $N$ and $k\ge1$ and reals $\xi>0$ and $0<p\le1$ be given. We say that a graph $G_r(N)$  has the property $\mathscr{C}_{N,p}^k(\xi)$ if for every $U\subseteq V^{(\ell+1)}$ and every family $\mathcal{F}_k\subseteq\mathcal{K}$ of pairwise disjoint $k$-sets with

\hspace*{.11cm}(i)\,\, $|\mathcal{F}_k|\le \xi N$ \,\,and

(ii)\,\, $|U|\le |\mathcal{F}_k|$,
\\
we have $\;e_{\Gamma}(\mathcal{F}_k,U)< 6p^k|U||\mathcal{F}_k|$.
\end{definition}

The following fact tells that for $G_r(N,p)$, {\bf a.a.s.} the corresponding graph $\Gamma(k,G_r(N,p))$  has no dense subgraph.
\begin{fact}\label{congestion}
For every integer $k\ge 1$ and real $\xi>0$, there exists $C>1$ such that if $p>C\big(\frac{\log N}{N}\big)^{1/k}$, then $\mathbb{P}(G_r(N,p)\in {\mathscr{C}}^k_{N,p}(\xi))=1-o(1)$.
\end{fact}

The following corollary follows immediately from Fact \ref{congestion}
since $\big(\frac{\log N}{N}\big)^{1/\Delta}\ge\big(\frac{\log N}{N}\big)^{1/k}$ for $1\le k\le \Delta$.
\begin{corollary}\label{congestion-2}
For every integer $\Delta\ge 1$ and real $\xi>0$, there exists $C>1$ such that if $p>C\big(\frac{\log N}{N}\big)^{1/\Delta}$, then $\mathbb{P}(G_r(N,p)\in \bigcap_{k=1}^\Delta{\mathscr{C}}^k_{N,p}(\xi))=1-o(1)$.
\end{corollary}\label{congestion-gene}
{\bf Proof of Fact \ref{congestion}.}
Let $\mathcal{F}_k$ and $U$ be defined as in Definition \ref{congestion-def}.
Note that if $X$ is a binomial random variable $Bi(M,q)$ then $\mathbb{P}(X\ge t)\le q^t{M\choose t}\le(eqM/t)^t$.
So for fixed pair $(\mathcal{F}_k,U)$, we have
\begin{align*}
\mathbb{P}\big(e_\Gamma(\mathcal{F}_k,U)\ge6p^k|U||\mathcal{F}_k|\big)
\le \Big(\frac{e}{6}\Big)^{6\xi Np^k|\mathcal{F}_k|}.
\end{align*}
Moreover, the number of choices for the pair $(\mathcal{F}_k,U)$ is at most
\[
\sum_{f=1}^{\xi N}\sum_{u=1}^{f}{{{r\choose k}N^k}\choose f}{N\choose u}.
\]
Therefore, the probability $G_r(N,p)\not\in {\mathscr{C}}^k_{N,p}(\xi)$ can be bounded from above by
\begin{align*}
\sum_{f=1}^{\xi N}\sum_{u=1}^{f}\exp&\Big(kf\log (rN)+u\log N-6\xi Np^kf\log (2\xi N/u)\Big)
\\&\le\sum_{f=1}^{\xi N}\sum_{u=1}^{f}\exp\Big((2k+1-6\xi C^k)f\log N\Big),
\end{align*}
which will tends to $0$ as $N\to\infty$ if we choose $C$ such that
$C^k>\frac{k+2}{3\xi}.$
\hfill$\Box$

\medskip
Now, let us tend to discuss the last property we need for the random $r$-partite graph $G_r(N,p)$. Let us first give the definitions for classes $\mathcal{B}_p^I$ and $\mathcal{B}_p^{II}$ of ``bad" tripartite graphs.
\begin{definition}\label{here-B}
 Let integers $m_1,m_2$ and $m_3$ and reals $\alpha,\epsilon',\epsilon,\mu>0$ and $0<p\le1$ be given.
\begin{itemize}
  \item
   { Let $\mathcal{B}_p^I(m_1,m_2,m_3,\alpha,\epsilon',\epsilon,\mu)$ be the family of tripartite graphs with three color classes $X,Y$ and $Z$, where $|X|=m_1$, $|Y|=m_2$ and $|Z|=m_3$, satisfying

     {\em($a$)} $(X,Y)$ and $(Y,Z)$ are $(\epsilon,\alpha,p)$-dense pairs and

{\em($b$)} there exists $X'\subseteq X$ with $|X'|\ge\mu|X|$ such that for every $x\in X'$ the pair $(N(x)\cap Y,Z)$ is not
$(\epsilon',\alpha,p)$-dense. }

  \item
   { Let $\mathcal{B}_p^{II}(m_1,m_2,m_3,\alpha,\epsilon',\epsilon,\mu)$ be the family of tripartite graphs with three color classes $X,Y$ and $Z$, where $|X|=m_1$, $|Y|=m_2$ and $|Z|=m_3$, satisfying

     {\em($a$)} $(X,Y)$, $(X,Z)$ and $(Y,Z)$ are $(\epsilon,\alpha,p)$-dense pairs and

{\em($b$)} there exists $X'\subseteq X$ with $|X'|\ge\mu|X|$ such that for every $x\in X'$ the pair $(N(x)\cap Y,N(x)\cap Z)$
is not $(\epsilon',\alpha,p)$-dense.
      }
\end{itemize}

\end{definition}

In the following, we define the family of graph $\mathscr{D}_{N,p}^\Delta$ that contains no element of $\mathcal{B}_p^I\cup \mathcal{B}_p^{II}$.
\begin{definition}\label{here-def}
 For integers $N$ and $\Delta\ge2$ and reals $\alpha,\gamma,\epsilon',\epsilon,\mu>0$ and $0<p\le1$ we say that a graph $G_r(N)$ has the property
 $\mathscr{D}_{N,p}^\Delta(\gamma,\alpha,\epsilon',\epsilon,\mu)$ if $G_r(N)$ contains no member from

\[
\mathcal{B}_p^I(m_1^I,m_2^I,m_3^I,\alpha,\epsilon',\epsilon,\mu)\cup\mathcal{B}_p^{II}(m_1^{II},m_2^{II},m_3^{II},\alpha,\epsilon',\epsilon,\mu)
\]
 with $m_1^I,m_3^I\ge \gamma p^{\Delta-1}N$ and $m_2^I,m_1^{II},m_2^{II},m_3^{II}\ge\gamma p^{\Delta-2}N$ as a subgraph.
\end{definition}

The following fact tells us that the random $r$-partite graph $G_r(N,p)$ has the property $\mathscr{D}_{N,p}^\Delta$ with high probability for suitable $p$.
\begin{fact}\label{here-p}
For every integer $\Delta\ge2$ and positive reals $\alpha,\epsilon'$ and $\mu$ there exists
$\epsilon=\epsilon(\Delta,\alpha,\epsilon',\mu)>0$
such that for every $\gamma>0$ there exists $C=C(\Delta,\alpha,\epsilon',\mu,\gamma)>1$ such that if $p>C(N/\log N)^{1/\Delta}$, then
$\mathbb{P}(G_r(N,p)\in {\mathscr{D}}^\Delta_{N,p}(\gamma,\alpha,\epsilon',\epsilon,\mu))=1-o(1)$.
\end{fact}

Before giving a proof for Fact \ref{here-p}, we have the following corollary.
\begin{corollary}\label{here-p-gen}
For all integers $\Delta,\bar{\Delta}\ge2$ and all reals $\alpha,\mu,\gamma$ and $\epsilon^*>0$, there exists  $C>1$ and $\epsilon_0,\epsilon_1,\dots,\epsilon_{\bar{\Delta}}$ satisfying $0<\epsilon_0\le\epsilon_1\le\dots\le\epsilon_{\bar{\Delta}}\le\epsilon^*$
such that  if $p>C(N/\log N)^{1/\Delta}$, then
$\mathbb{P}(G_r(N,p)\in \bigcap_{k=1}^{\bar{\Delta}}{\mathscr{D}}^\Delta_{N,p}(\gamma,\alpha,\epsilon_k,\epsilon_{k-1},\mu))=1-o(1)$.
\end{corollary}
{\bf Proof.} Let $\Delta,\bar{\Delta}\ge2$ and all reals $\alpha,\mu,\gamma$ and $\epsilon^*>0$ be given. First, set $\epsilon_{\bar{\Delta}}=\epsilon^*$. For $\bar{\Delta}>k\ge1$,  set recursively that $\epsilon_{k-1}=\min\{\epsilon_k,\epsilon(\Delta,\alpha,\epsilon_k,\mu)\}$, where $\epsilon(\Delta,\alpha,\epsilon_k,\mu)$ is given by Fact \ref{here-p}. Then set $C$ to be the maximum of all $C(\Delta,\alpha,\epsilon_k,\mu,\gamma)$ for $k=1,\dots,\bar{\Delta}$.
Thus, Fact \ref{here-p} guarantees that {\bf a.a.s.} $G_r(N,p)\in {\mathscr{D}}^\Delta_{N,p}(\gamma,\alpha,\epsilon_k,\epsilon_{k-1},\mu)$ for $1\le k\le\bar{\Delta}$.\hfill$\Box$

\medskip
We first consider Fact \ref{here-p} for the special case. Define
$$\mathcal{B}_p^I(m,\alpha,\epsilon',\epsilon,\mu)=\mathcal{B}_p^I(pm,m,pm,\alpha,\epsilon',\epsilon,\mu) $$
and
$$
\mathcal{B}_p^{II}(m,\alpha,\epsilon',\epsilon,\mu)=\mathcal{B}_p^{II}(m,m,m,\alpha,\epsilon',\epsilon,\mu)
$$
Similarly, for integers $N$ and $\Delta\ge2$ and reals $\alpha,\gamma,\epsilon',\epsilon,\mu>0$ and $0<p\le1$, we say that a graph $G_r(N)$ has the property
 $\widehat{\mathscr{D}}_{N,p}^\Delta(\gamma,\alpha,\epsilon',\epsilon,\mu)$ if $G_r(N)$ contains no member from $\mathcal{B}_p^I(m,\alpha,\epsilon',\epsilon,\mu)\cup\mathcal{B}_p^{II}(m,\alpha,\epsilon',\epsilon,\mu)$ as a subgraph.

\begin{lemma}\label{here-p-pt}
For $\Delta\ge2$ and  $\alpha,\epsilon',\mu \in (0,1]$ there exists $\epsilon>0$
such that for every $\gamma\in(0,1]$ there exists $C\ge1$ such that if $p>C(N/\log N)^{1/\Delta}$, then
$\mathbb{P}(G_r(N,p)\in {\widehat{\mathscr{D}}}^\Delta_{N,p}(\gamma,\alpha,\epsilon',\epsilon,\mu))=1-o(1)$.
\end{lemma}
{\bf Proof.} Let $\Delta\ge2$ and  $\alpha,\epsilon',\mu\in (0,1]$ be given, and let
\begin{align}\label{beta}
\beta=\frac{\alpha^2}{4e^2}\Big(\frac{1}{e}\Big)^{8/(\alpha\mu)}.
\end{align}
For $\alpha,\beta$ and $\epsilon'$, there exist $\epsilon_1$ and $L_1$, and $\epsilon_2$ and $L_2$ according to Theorem \ref{num-X'} and Corollary \ref{num-den-pair}, respectively.
 Fix
$
\epsilon=\min\{\alpha/2,\mu/4,\epsilon_1,\epsilon_2\},
$
and for every $\gamma>0$ we set
$
C=({4}/{\gamma})^{1/\Delta},
$
and let $N$ be sufficiently large.

\smallskip
First we show that {\bf a.a.s.} $G_r(N,p)$ contains no element from
\[
\mathcal{B}_p^{II}(m,\alpha,\epsilon',\epsilon,\mu).
\]

Let $T$ be the tripartite graph with color classes $X,Y$ and $Z$ from $\mathcal{B}_p^{II}(m,\alpha,\epsilon',\epsilon,\mu)$.
From Definition \ref{here-B}, the bipartite subgraphs $T[X,Y]$, $T[X,Z]$ and $T[Y,Z]$ of $T$ contain at least
$(\alpha-\epsilon)pm^2$ edges each. Also, there is a set $X'\subseteq X$ with $|X'|\ge\mu|X|$ such
that for every $x\in X'$ the pair $(N_T(x)\cap Y,N_T(x)\cap Z)$ is not $(\epsilon',\alpha,p)$-dense. Set
\[
X''=\{x\in X': |N_T(x)\cap Y|\ge\alpha pm/2\,\,\mbox{and}\,\,|N_T(x)\cap Z|\ge\alpha pm/2\}.
\]
Since $(X,Y)$ and $(X,Z)$ are $(\epsilon,\alpha,p)$-dense, we have
$
|X''|\ge(1-2\epsilon/\mu)|X'|\ge|X'|/2\ge\mu pm/2.
$
Choose a subset $X'''\subseteq X''$ with  $|X'''|=\mu pm/2$.

Fix $x\in X'''$. Similarly, there are sets $Y_x'\subseteq N_T(x)\cap Y$ and $Z_x'\subseteq N_T(x)\cap Z$
of size $\epsilon'\alpha pm/2$ such that $d_{T,p}(Y_x',Z_x')<\alpha-\epsilon'$. Now let $Y_x$ and $Z_x$ be such
that $Y_x'\subseteq Y_x\subseteq N_T(x)\cap Y$ and $Z_x'\subseteq Z_x\subseteq N_T(x)\cap Z$ with $|Y_x|=|Z_x|=\alpha pm/2$. Then,  $(Y_x, Z_x)$ is not
$(\epsilon',\alpha,p)$-dense. Thus, we have a family $\{(Y_x, Z_x): x\in X'''\}$ such that each pair in which is not $(\epsilon',\alpha,p)$-dense.

Note that $\alpha pm/2\ge L_2/p$ if $N$ is large. Hence, apply Corollary \ref{num-den-pair} to $T[Y,Z]$, the probability that $T[X''',Y,Z]$ appears in $G_r(N,p)$ can be bounded from above by
\begin{align*}
\sum_{t\ge(\alpha-\epsilon)pm^2}&{r\choose 3}{N\choose m}^3{m^2\choose t}p^t\cdot\left(\beta^{\alpha pm/2}{m\choose \alpha pm/2}^2\right)^{\mu pm/2}p^{\mu\alpha pm^2/2}
\\&\le \sum_{t\ge(\alpha-\epsilon)pm^2}{r\choose 3}{N\choose m}^3\left(\frac{pm^2e}{t}\right)^t\cdot\left(\beta^{1/2}\frac{e}{\alpha/2}\right)^{\mu\alpha pm^2/2}
\\&<pm^2r^3N^{3m}\left(e\left(\frac{2e}{\alpha}\right)^{\mu\alpha/2}\beta^{\mu\alpha/4}\right)^{pm^2},
\end{align*}
where the last inequality holds as the function $f(t)=(pm^2e/t)^t$ is maximized for $t=pm^2$.
Note that $e(\frac{2e}{\alpha})^{\mu\alpha/2}\beta^{\mu\alpha/4}= 1/e$ and $pm^2$ is much larger than $ m\log N$, we have the right-hand side of the last inequality tends to $0$ as $N\rightarrow\infty$.

\smallskip
The fact that {\bf a.a.s.} $G_r(N,p)$ contains no element from
$\mathcal{B}_p^{I}(m,\alpha,\epsilon',\epsilon,\mu)$ is similar (indeed a little simpler) as above except some differences in calculation, and so we omit the proof.
 \hfill$\Box$

\medskip
Now, Fact \ref{here-p} follows immediately from Lemma \ref{here-p-pt} and the following result obtained in (\cite{Kohayakawa}, Claim 18 and the remark afterwards).
\begin{lemma}
For an integer $\Delta\ge 2$ and positive reals $\alpha,\epsilon',\mu$ and $\widehat{\epsilon}$ there exists $\epsilon>0$ such that for every
$\gamma>0$ there exists $C > 1$ and $N_0$ such that if
$N\ge N_0$ and $p > C(\log N/N)^{1/\Delta}$, then every tripartite graph $T\in \mathcal{B}_p^I(m_1^I,m_2^I,m_3^I,\alpha,\epsilon',\epsilon,\mu)\cup\mathcal{B}_p^{II}(m_1^{II},m_2^{II},m_3^{II},\alpha,\epsilon',\epsilon,\mu)$ with $m_1^I,m_3^I\ge \gamma p^{\Delta-1}N$ and $m_2^I,m_1^{II},m_2^{II},m_3^{II}\ge\gamma p^{\Delta-2}N$ contains a subgraph $\widehat{T}\in\mathcal{B}_p^I(m,\alpha,\epsilon',\epsilon,\mu)\cup\mathcal{B}_p^{II}(m,\alpha,\epsilon',\epsilon,\mu)$.
\end{lemma}

\section{Sparse partition universal multipartite Graphs for $\mathcal{H}(\Delta,n)$}
In this section, we will show Theorem \ref{main}, i.e., for fixed $\Delta\ge2$, we can choose suitable constants $r=r(\Delta)$, $B=B(\Delta)$ and $C=C(\Delta)$ such that for any $H\in\mathcal{H}(\Delta,n)$ if $N\ge Bn$ and $p=C(\log N/N)^{1/\Delta}$, then {\bf a.a.s.} the random graph $G_r(N,p)\rightarrow (H,H)$.  A key tool for the proof is an elegant embedding method developed in \cite{Kohayakawa}.
\begin{lemma}{}\label{main-le}
For every $\Delta\ge2$ there exist $\bar{\Delta}\ge2$ and positive constants $\mu,\alpha,\epsilon^*,\xi$ and $\gamma>0$ and $B>1$ and $n_0$ such that for every $\epsilon_0,\dots,\epsilon_{\bar{\Delta}}$ satisfying $0<\epsilon_0\le\dots\le\epsilon_{\bar{\Delta}}\le\epsilon^*$ and for every $n\ge n_0$ the following holds.
If $G_r(N)$ is a $r$-partite graph has $r$-color classes $V^{(1)},V^{(2)},\dots,V^{(r)}$ with $|V^{(i)}|=N\ge Bn$ ($1\le i\le r$) such that for some $0< p\le1$ we have

\hspace*{.2cm}{(i)} $G_r(N)\in{\mathscr{U}}_{N,p}$,

\hspace*{.1cm}{(ii)} $G_r(N)\in{\mathscr{C}}^k_{N,p}(\xi)$ for every $k=1,\dots,{\Delta}$, and

{(iii)} $G_r(N)\in {\mathscr{D}}^\Delta_{N,p}(\gamma,\alpha,\epsilon_k,\epsilon_{k-1},\mu)$ for every $k=1,\dots,\bar{\Delta}$,

\smallskip
\noindent
then $G_r(N)\rightarrow (H,H)$ for $H\in\mathcal{H}(\Delta,n)$.
\end{lemma}

Before we prove Lemma \ref{main-le}, we deduce Theorem \ref{main} from it.

\smallskip
\noindent{\bf Proof of Theorem \ref{main}.} Let $\Delta\ge2$ be given, Lemma \ref{main-le}
yields constants $\bar{\Delta}\ge2$ and $\mu,\alpha,\epsilon^*,\xi$ and $\gamma>0$ and $B>1$ and $n_0$.
Moreover, from Fact \ref{uniform}, Corollary \ref{congestion-2} and Corollary \ref{here-p-gen},
there exists a constant $C$ such that for $p>C(\log N/N)^{1/\Delta}$ the
random $r$-partite graph $G_r(N,p)$ {\bf a.a.s.} satisfies the assumptions ($i$), ($ii$) and ($iii$).
Therefore, Lemma \ref{main-le} asserts that {\bf a.a.s.} $G_r(N,p)\to (H,H)$  for every $H\in\mathcal{H}(\Delta,n)$
whenever $N\ge Bn$, which completes the proof of Theorem \ref{main}.\hfill$\Box$

\smallskip
In order to prove Lemma \ref{main-le}, we also need the following result which relates to Tur\'{a}n number for $K_r$ in complete $r$-partite graph $K_r(k)$.
Although the proof of the following lemma is simple, it is crucial in the procedure of embedding bounded degree graph to random multipartite graph.
\begin{lemma}\label{turan}
For integers $k\ge 1$ and $r\ge 2$, let $t_r(k)$ be the maximum number of edges in a subgraph of $K_r(k)$ that contains no $K_r$. We have
\[
t_r(k)= \left[{r\choose 2}-1\right] k^2.
\]
\end{lemma}
{\bf Proof.} The lower bound for $t_r(k)$ follows by deleting all edges between a pair of color classes of $K_r(k)$.
On the other hand, we shall prove  by induction of $k$ that if a subgraph $G=G(V^{(1)},\dots,V^{(r)})$ of $K_r(k)$ contains no $K_r$, then $e(G)\le \left[{r\choose 2}-1\right] k^2$. Suppose $k\ge 2$ and $r\ge3$ as it is trivial for $k=1$ or $r=2$. Now, suppose  that
$G$ has the maximum possible number of edges subject to this condition. Then $G$ must contain $K_{r}-e$ as a
subgraph, otherwise we could add an edge and the resulting graph would still not contain $K_r$. Denote the vertex set of this $K_{r}-e$ by $X$, we have $|X\cap V^{(i)}|=1$ for $i=1,2,\dots,r$. Without loss of generality, suppose $e=\{v_1,v_2\}$, where $v_1\in V^{(1)}$ and $v_2\in V^{(2)}$. Note that $G$ contains no $K_r$, we have for $i=1,2$ there is no vertex in $V^{(i)}\setminus\{v_i\}$ is adjacent to all the vertices of $X\setminus\{v_i\}$. Thus, together with the induction hypothesis, we can deduce that there are at least
\[
(k-1)^2+2(k-1)+1=k^2
\]
edges should be deleted from $K_r(k)$, which completes the induction step hence the proof.  \hfill  $\Box$

\smallskip
Now, we are ready to give a proof for Lemma \ref{main-le}.

\smallskip
\noindent{\bf Proof of Lemma \ref{main-le}.} The proof consists of four parts. Firstly, we fix all
constants needed in the proof. Secondly, we consider the given $r$-partite graph $G_r(N)$
along with a fixed 2-coloring of its edges. In order to embed every
graph $H\in\mathcal{H}(\Delta,n)$ into one of the two monochromatic subgraphs of $G_r(N)$,
we first prepare the graph $G_r(N)$ and here the sparse multipartite regularity lemma will be the key
tool. Thirdly, we shall prepare a given graph $H\in\mathcal{H}(\Delta,n)$ for the embedding.
Finally, we will embed $H$ into a monochromatic subgraph of $G_r(N)$.

Let $\Delta\ge2$ be fixed, and let $\bar{\Delta}=\Delta^4+\Delta$ and $r=R(K_{\bar{\Delta}})$.
The constants $\mu,\alpha,\epsilon^*,\xi,\gamma,B$ and $n_0$ involved in the proof of Lemma \ref{main-le} are defined as follows. Set
\begin{align}\label{mu-alph-epi*}
\mu=\frac{1}{4\Delta^2},\,\,\,\alpha=\frac{1}{3},\,\,\,\mbox{and}\,\,\,\epsilon^*=\frac{1}{6{\Delta}}.
\end{align}
Let
\begin{align}\label{epi-K-t0}
\epsilon=\frac{\epsilon_0}{2},\,\,\,\lambda=2,\,\,\,\mbox{and}\,\,\, t_0=1
\end{align}
be given, there exist $T_0,\eta$ and $N_0$ that are guaranteed by the sparse multipartite regularity lemma,
Theorem \ref{sparse-reg}. Finally, we set
\begin{align}\label{xi-ga-B}
\xi=\frac{1}{6\cdot4^{\Delta+1}\cdot T_0},\,\,\,\gamma=\frac{1-\epsilon}{4^{\Delta-1}T_0},\,\,\, B=\frac{1}{\xi},
\end{align}
and
\begin{align}\label{n0}
n_0=\max\Big\{\frac{N_0}{B}, \,\, e^{1/\eta}, \,\, e^{{T_0}/{\epsilon(1-\epsilon)}}\Big\}.
\end{align}

\medskip
Now, let $\epsilon_0,\dots,\epsilon_{\bar{\Delta}}$ satisfy
\begin{align}\label{eps0...}
0<\epsilon_0\le\dots\le\epsilon_{\bar{\Delta}}\le\epsilon^*=\frac{1}{6{\Delta}}
\end{align}
and let $n\ge n_0$ be given. Let $G_r(N)$ be the $r$-partite graph which has $r$-color classes $V^{(1)},\dots,V^{(r)}$ with $|V^{(i)}|=N$ ($1\le i\le r$),
where $N\ge Bn\ge N_0$, satisfies assumptions ($i$)-($iii$) of Lemma \ref{main-le}. Denote $V=\bigcup_{i=1}^r V^{(i)}$. Consider an
arbitrary red/blue edge coloring of $G_r(N)$, and let $G_R = (V,E_R)$ and $G_B = (V,E_B)$ be the corresponding monochromatic
subgraphs.

In the following, we apply the sparse multipartite regularity lemma with $\epsilon=\epsilon_0/2$,
$\lambda = 2$, $t_0 = 1$, and some $p$ to $G_R$.
From property ($i$) of Lemma \ref{main-le}, the graph $G_r(N)$ is $(1/\log N, 1+1/\log N)$-bounded (see Definition \ref{e(UiUj)}). Since $G_R\subseteq G_r(N)$ and
$1/\log N\le\eta\le1$ from (\ref{n0}), we have $G_R$ is $(\eta,\lambda)$-bounded.

Consequently, from Theorem \ref{sparse-reg}, we have an partition
$\big\{V^{(i)}_0,V^{(i)}_1,\cdots, V^{(i)}_t\big\}$ for each $V^{(i)}$, where $t$ is same for each $V^{(i)}$ and $t_0\le t \le T_0$, such that

\vspace*{-.5cm}
\begin{align*}
(i)\,\, &|V^{(i)}_0|\le\epsilon N\,\,\mbox{and}\,\,|V^{(i)}_1|=\dots=|V^{(i)}_t|\,\,\mbox{for}\,\, 1\le i\le r;
\\(ii)\,\, & \mbox{All but at most}\,\,\epsilon t^2\,\,\mbox{pairs}\,\,\big(V^{(i)}_a,V^{(j)}_b\big)\,\,\mbox{are}\,\, (\epsilon,p)-regular\,\,\mbox{for}\,\,1\le i,j\le r\,\,\mbox{and}\,\, 1\le a,b \le t.
\end{align*}
Let $F$ be the subgraph of $r$-partite complete graph $K_r(t)$, whose vertices are $\big\{V^{(i)}_a\;|\;1\le a\le t, 1\le i\le r\big\}$ in which a pair  $\big(V^{(i)}_a,V^{(j)}_b\big)$ for $i\neq j$ is adjacent if and only if
the pair  is $(\epsilon,p)$-regular in $G_R$.  Then the number of edges of
$F$ is at least
\[
t^2{r\choose 2}- \epsilon t^2>\left[{r\choose 2}-1\right]t^2=t_r(t).
\]
Hence, Lemma \ref{turan} implies that $F$ contains a complete graph $K_r$ with $r$ vertices.
 Without loss of generality, assume that $V^{(1)}_1,\dots, V^{(r)}_1$ are pairwise $(\epsilon,p)$-regular for $G_R$.
 For convenience, let us denote $V^{(1)}_1,\dots, V^{(r)}_1$ by $A_1,\dots, A_r$.
Moreover, since $G_r(N)\in{\mathscr{U}}_{N,p}$, we have $|d_{G_R,p}(A_i,A_j)+ d_{G_B,p}(A_i,A_j)-1|\le1/\log N$.
Note from (\ref{n0}) that $N/\log N\le\epsilon(1-\epsilon) N/T_0\le\epsilon |A_i|$, we have $(A_i,A_j)$ is $(\epsilon+ 2/\log N, p)$-regular for
the graph $G_B$. From (\ref{epi-K-t0}) and (\ref{n0}), we have $\epsilon+ 2/\log N\le \epsilon_0$ and,
hence, $(A_i,A_j)$ is $(\epsilon_0,p)$-regular both for $G_R$ and $G_B$ for all $\{i,j\}\in{[r]\choose2}$. Therefore,
\[
\max\{d_{G_R,p}(A_i,A_j),d_{G_B,p}(A_i,A_j)\}\ge\frac{1}{2}-\frac{1}{2\log N}\ge\frac{1}{3}
\]
for every $\{i,j\}\in{[r]\choose2}$.

Color an edge $\{i,j\}\in{[r]\choose2}$ red if $d_{G_R,p}(A_i,A_j)\ge d_{G_B,p}(A_i,A_j)$ and blue otherwise.
Since $r\ge R(K_{\bar{\Delta}})$, we can suppose without loss of generality that there exists a monochromatic red clique
$K_{\bar{\Delta}}$ with vertex set $[\bar{\Delta}]\subseteq[r]$ satisfying
\begin{align}\label{dense>1/3}
(A_i,A_j)\,\,\mbox{ is}\,\, (\epsilon_0,p)\mbox{-regular for}\,\, G_R \,\,\mbox{and}\,\, d_{G_R,p}(A_i,A_j)\ge 1/3 \,\,\mbox{for all}\,\, \{i,j\}\in{[\bar{\Delta}]\choose2},
\end{align}
 and we shall show that $G_R$ induced on $\bigcup_{i\in [\bar{\Delta}]}A_i$ contain any $H\in\mathcal{H}(\Delta,n)$.

\medskip
\noindent{\bf Preparing $H$.} \, Fix $H=(W,E(H))\in\mathcal{H}(\Delta,n)$. Define the third power $H^3=(W,E(H^3))$ of $H$ such that
$\{w,w'\}\in E(H^3)$ if and only if $w\neq w'$ and there exists a $w-w'$-path with at most three edges in $H$. Since $\Delta(H)\le\Delta$, we have
$\Delta(H^3)\le\Delta+\Delta(\Delta-1)+\Delta(\Delta-1)^2=\Delta^3-\Delta^2+\Delta$. Hence, $H^3$ can not be a complete graph as $|W|=n$ is large.
Consequently, $\chi(H^3)\le\Delta(H^3)\le\Delta^3-\Delta^2+\Delta$. Let $``f"$ be a $(\Delta^3-\Delta^2+\Delta)$-vertex
coloring of $H^3$.
We say two vertices $w$ and $w'$ are equivalent according to $``f"$ if $f(w) = f(w')$, and
\[
|x\in N_H(w):\,f(x)<f(w)|=|x\in N_H(w):\,f(x)<f(w')|,
\]
i.e., with same ``left-degrees" --- the
number of neighbors with colors of smaller number. This equivalence relation partitions $W$ into at most $(\Delta^3-\Delta^2+\Delta)(\Delta+1)=\bar{\Delta}$ classes
as each vertex has ``left-degrees" at most $\Delta$. Denote the partition classes by $W_1,\dots,W_{\bar{\Delta}}$
(may have empty classes) and let $g: W\rightarrow [\bar{\Delta}]$ be the corresponding partition function, i.e.,
$g(w)=j\,\,\mbox{if and only if}\,\,w\in W_j.$
Thus, if $w,w'\in W_j$, then
$
|x\in N_H(w):\,g(x)<g(w)|=|x\in N_H(w):\,g(x)<g(w')|.
$
For $\ell\le g(w)$, denote the ``{\em left-degree}" of $w$ { with respect to} $g$ { and} $\ell$ by
\[
\mbox{ldeg}_g^\ell(w)=|x\in N_H(w):\,g(x)\le\ell|.
\]
Note that if $w,w'\in W_j$, then their distance in $H$ is at least four, which implies
\begin{align}\label{NwNw'=0}
|N_H(w)\cap N_H(w')|=0\,\,\mbox{and}\,\,e_H(N_H(w), N_H(w'))=0.
\end{align}

\medskip
\noindent{\bf Embedding of $H$ to $G_R$.}\, Now,  we will embed the vertex class $W_\ell$ into $A_\ell$
one at a time, for $\ell=1,\dots,\bar{\Delta}$. Indeed, we shall inductively  verify the following statement ($\mathcal{S}_\ell$)
for $\ell=0,1,\dots,\bar{\Delta}$.

\smallskip
($\mathcal{S}_\ell$)  There is an embedding $\varphi_\ell$ of $H[\cup_{j=1}^\ell W_j]$ into $G_R[\cup_{j=1}^\ell A_j]$
such that for every $z\in \cup_{j=\ell+1}^{\bar{\Delta}}W_j$,
there exists a {\em candidate set} $C_\ell(z)\subseteq A_{g(z)}$ given by

\[
(a)\,\, C_\ell(z)=\Bigg(\bigcap_{x\in N_H(z),\, g(x)\le \ell}N_{G_R}(\varphi_\ell(x))\Bigg)\cap A_{g(z)},
\]
satisfying

($b$) $|C_\ell(z)|\ge(p/4)^{\mbox{ldeg}_g^\ell(z)}m$, where $m=|A_{g(z)}|\ge(1-\epsilon)N/t$, and

($c$) $\big(C_\ell(z),C_\ell(z')\big)$ is $(\epsilon_\ell,1/3,p)$-dense in $G_R$ for any edge $\{z,z'\}\in E(H)$ with $g(z),g(z')>\ell$.

\smallskip
\noindent{\em Remark.} The definition of $C_\ell(z)$ implies that if we embed $z$ into $C_\ell(z)$, then its image will be adjacent to all vertices $\varphi_\ell(x)$ with $x\in N_H(z)\cap(W_1\cup\dots\cup W_\ell)$.

Note that ($\mathcal{S}_{\bar{\Delta}}$) gives an embedding $\varphi_{\bar{\Delta}}$ of $H[\cup_{j=1}^{\bar{\Delta}} W_j]$ into $G_R[\cup_{j=1}^{\bar{\Delta}} A_j]$,
so $G_R$ contains $H$ as a subgraph. Thus, verifying ($\mathcal{S}_{\ell}$) inductively for $\ell=0,1,\dots,\bar{\Delta}$ completes the proof of Lemma \ref{main-le}. In the following, we shall use the letter ``$x$" for vertices that
have been embedded, ``$y$" for that will be embedded in the
current step, while ``$z$" for that we shall embed at a later step.

Let us verify ($\mathcal{S}_0$) at first. In this case, $\varphi_0$ is the empty
mapping and for every $z\in W$, we have $C_0(z) = A_{g(z)}$ according to ($a$) as there is no vertex $x\in N_H(z)$ with $g(x)\le 0$.
So property ($b$) follows directly, and property ($c$) follows from (\ref{dense>1/3}).

For the inductive step, suppose ($\mathcal{S}_\ell$) holds for $\ell<\bar{\Delta}$,
and we shall verify ($\mathcal{S}_{\ell+1}$) by embedding $W_{\ell+1}$ into $A_{\ell+1}$ with the required properties.
Note from (\ref{NwNw'=0}) that $|N_H(z)\cap W_{\ell+1}|\le1$ for every $z\in \cup_{j=\ell+2}^{\bar{\Delta}}W_j$.
Hence, for every ``right-neighbor" $z\in\bigcup_{j=\ell+2}^{\bar\Delta}W_j$ of $y\in W_{\ell+1}$, the new candidate set will be $C_{\ell+1}(z)=C_\ell(z)\cap N_{G_R}(\varphi_{\ell+1}(y))$.
For each $y \in W_{\ell+1}$, we should find
a suitable subset $C(y) \subseteq C_\ell(y)$ such that if $\varphi_{\ell+1}(y)$ is
chosen from $C(y)$, then the new candidate set $C_{\ell+1}(z)$ will satisfy properties ($b$) and ($c$)
of ($\mathcal{S}_{\ell+1}$).
However, since in general $|C(y)| \le |C_\ell(y)|=o(|A_{g(z)}|)$ if $\mbox{ldeg}_g^\ell(y)\ge 1$,
 we should select $\varphi_{\ell+1}(y)$ from $C(y)$ carefully such that if $y\neq y'$ then $\varphi_{\ell+1}(y)\neq \varphi_{\ell+1}(y')$.
Here we shall apply Hall's condition, a similar idea was used in \cite{alon-furedi,rodl-ruci,Kohayakawa}. The details are contained in the following two claims.

Fix $y \in W_{\ell+1}$ and denote the ``right neighbors" of $y$ by
\[
N_H^{\ell+1}(y)=\{z\in N_H(y):\,g(z)>\ell+1\}.
\]
For a vertex $v\in C_\ell(y)$, let $\widehat{C}_\ell(z)=N_{G_R}(v)\cap C_\ell(z)$ if $z\in N_H^{\ell+1}(y)$ and $\widehat{C}_\ell(z)=C_\ell(z)$ if $z\not\in N_H^{\ell+1}(y)$.

\medskip
\noindent{\bf Claim 1.} For every $y \in W_{\ell+1}$, there exists a subset $C(y)\subseteq C_\ell(y)$ with
$|C(y)|\ge (1-\Delta\epsilon_{\ell}-\Delta^2\mu)|C_\ell(y)|$
 such that for every $v\in C(y)$, the following ($b'$) and ($c'$) hold.

($b'$) $|N_{G_R}(v)\cap C_\ell(z)|\ge({p}/{4})^{\mbox{ldeg}_g^{\ell+1}(z)}|A_{g(z)}|$ for every $z\in N_H^{\ell+1}(y)$,

($c'$) $(\widehat{C}_\ell(z),\widehat{C}_\ell(z'))$ is $(\epsilon_{\ell+1}, 1/3, p)$-dense for all edges $\{z,z'\}$ of $H$ with $g(z)$, $g(z')>\ell+1$ and
$\{z,z'\}\cap N_H^{\ell+1}(y)\neq\emptyset$.

\medskip
\noindent{\bf Proof.}
Fix $z\in N_H^{\ell+1}(y)$. Since $(C_\ell(y),C_\ell(z))$
is an $(\epsilon_\ell, 1/3, p)$-dense pair from ($c$) of ($\mathcal{S}_{\ell}$), we have there exist at most $\epsilon_\ell|C_\ell(y)|$ vertices $v\in C_\ell(y)$ such that
\[
|N_{G_R}(v)\cap C_\ell(z)|<\Big(\frac{1}{3}-\epsilon_{\ell}\Big)p|C_\ell(y)|.
\]
Note that $|N_H^{\ell+1}(y)|\le\Delta$, we have from ($b$) and ($c$) of ($\mathcal{S}_{\ell}$)
that all  but at most $\Delta\epsilon_{\ell}|C_\ell(y)|$ vertices $v\in C_\ell(y)$   satisfy that  for every $z\in N_H^{\ell+1}(y)$,
\begin{align}\label{NcapC}
|N_{G_R}(v)\cap C_\ell(z)|\nonumber{\ge}\Big(\frac{1}{3}-\epsilon_{\ell}\Big)p\Big(\frac{p}{4}\Big)^{\mbox{ldeg}_g^{\ell}(z)}|A_{g(z)}|
\stackrel{\mathrm{(\ref{eps0...})}}{\ge}\Big(\frac{p}{4}\Big)^{\mbox{ldeg}_g^{\ell+1}(z)}|A_{g(z)}|.
\end{align}

Now, fix an edge $e = \{z, z'\}$ with $g(z), g(z') > \ell+1$ and
with at least one end vertex in $N_H^{\ell+1}(y)$. Clearly, there are at most $\Delta(\Delta-1)<\Delta^2$ such
edges.

If both vertices $z$ and $z'$ are neighbors of $y$, i.e., $z, z'\in N_H^{\ell+1}(y)$,
then
\[
\max\{\mbox{ldeg}_g^{\ell}(y), \mbox{ldeg}_g^{\ell}(z), \mbox{ldeg}_g^{\ell}(z')\}\le \Delta-2,
\]
since $y, z$, and $z'$ have at least two neighbors in $W_{\ell+1}\cup\dots\cup W_{\bar{\Delta}}$.
Hence, property ($b$) of ($\mathcal{S}_{\ell}$) implies
\[
\min\{|C_\ell(y)|, |C_\ell(z)|, |C_\ell(z')|\}\ge\Big(\frac{p}{4}\Big)^{ \Delta-2}(1-\epsilon)\frac{N}{T_0}\ge \gamma p^{\Delta-2}N.
\]
Recall that $\alpha= 1/3$, see (\ref{mu-alph-epi*}). Hence $G_R \subseteq G_r(N)$ and $G_r(N)\in {\mathscr{D}}^\Delta_{N,p}(\gamma,\alpha,\epsilon_{\ell+1},\epsilon_\ell,\mu)$
imply that there are at most $\mu|C_\ell(y)|$ vertices $v\in C_\ell(y)$ such that the
pair $(N_{G_R}(v)\cap C_\ell(z),N_{G_R}(v)\cap C_\ell(z'))$ is not $(\epsilon_{\ell+1}, 1/3, p)$-dense.

If, on the other hand, say, only $z\in N_H^{\ell+1}(y)$ and $z'\not\in N_H^{\ell+1}(y)$, then
\[
\max\{\mbox{ldeg}_g^{\ell}(y),\mbox{ldeg}_g^{\ell}(z')\}\le \Delta-1\,\,\mbox{and}\,\,\mbox{ldeg}_g^{\ell}(z)\le \Delta-2.
\]
Hence, similarly, we have
$\min\{|C_\ell(y)|,|C_\ell(z')|\}\ge \gamma p^{\Delta-1}N \,\,\mbox{and}\,\,  |C_\ell(z)|\ge \gamma p^{\Delta-2}N$.
Thus, the fact $G_r(N)\in {\mathscr{D}}^\Delta_{N,p}(\gamma,\alpha,\epsilon_{\ell+1},\epsilon_\ell,\mu)$
implies that there are at most $\mu|C_\ell(y)|$ vertices $v\in C_\ell(y)$ such that the
pair $(N_{G_R}(v)\cap C_\ell(z), C_\ell(z'))$ is not $(\epsilon_{\ell+1}, 1/3, p)$-dense.

Therefore, for either case, deleting all ``bad" vertices from $C_\ell(y)$, we find a subset $C(y)\subseteq C_\ell(y)$ with size
$|C(y)|\ge (1-\Delta\epsilon_{\ell}-\Delta^2\mu)|C_\ell(y)|$
 such that for every $v\in C(y)$, the following ($b'$) and ($c'$) hold.\hfill$\Box$

\smallskip

Now, let us turn to the second part of the inductive step. We shall choose $\varphi_{\ell+1}(y)\in C(y)$ such that
$\varphi_{\ell+1}(y)\neq \varphi_{\ell+1}(y')$ for two vertices $y, y'\in W_{\ell+1}$.
This can be achieved from the following claim.

\medskip
\noindent{\bf Claim 2.} There is an injective mapping
$
\psi: W_{\ell+1}\rightarrow \bigcup_{y\in W_{\ell+1}}C(y)
$
such that $\psi(y)\in C(y)$ for every $y\in W_{\ell+1}$.

\medskip
\noindent{\bf Proof.} It suffices to verify Hall's condition that for every $Y \subseteq W_{\ell+1}$,
\begin{align}\label{|Y|}
|Y|\le \Big|\bigcup_{y\in Y}C(y)\Big|.
\end{align}
Recall that $\mbox{ldeg}_g^{\ell}(y)=\mbox{ldeg}_g^{\ell}(y')$ for any two distinct vertices $y, y'\in W_{\ell+1}$ and so we can set
\begin{align}\label{k=ldeg-y}
k=\mbox{ldeg}_g^{\ell}(y) \,\,\mbox{for all}\,\, y\in W_{\ell+1}.
\end{align}
From Claim 1, property ($b$) of ($\mathcal{S}_\ell$), and (\ref{mu-alph-epi*}) and (\ref{eps0...}), we have that
\begin{align}\label{Cy}
|C(y)|\ge(1-\Delta\epsilon_{\ell}-\Delta^2\mu)|C_\ell(y)|\ge(1-\Delta\epsilon_{\ell}-\Delta^2\mu)\Big(\frac{p}{4}\Big)^k(1-\epsilon)\frac{N}{T_0}
\ge\frac{1}{4^{k+1}}p^k\frac{N}{T_0}.
\end{align}
Hence,
the assertion (\ref{|Y|}) holds if $|Y|\le {4^{-k-1}}p^k{N}/{T_0}$.
Thus, we suppose that $|Y|> {4^{-k-1}}p^k{N}/{T_0}$.
Denote the $k$-tuple $K(y)=\{u_1(y),\dots,u_k(y)\}$ by the neighbors of $y$ that have been embedded already. Clearly, $K(y)=N_H(y)\setminus N_H^{\ell+1}(y)$.
From (\ref{NwNw'=0}), we have  that
\begin{align*}
&|K(y)\cap W_i|\le1, \,\,\,\mbox{for}\,\,\,1\le i\le \ell;
\\K(y)\cap K(&y')=\emptyset, \,\,\,\mbox{for distinct vertices}\,\,\,y, y'\in W_{\ell+1}.
\end{align*}
Thus, the
sets of already embedded vertices $\varphi_\ell(K(y))$ and $\varphi_\ell(K(y'))$ are disjoint as well and,
therefore, $\mathcal{F}_k=\{\varphi_\ell(K(y)):\,y\in Y\}\subseteq {{V_1\cup\dots\cup V_\ell}\choose k}$  is a family of pairwise disjoint $k$-sets
with $|\varphi_\ell(K(y))\cap A_i|\le1$ for  $1\le i\le \ell$.
Note that
\[
C(y)\subseteq C_\ell(y)=\Big(\bigcap_{x\in K(y)}N_{G_R}(\varphi_\ell(x))\Big)\cap A_{\ell+1}.
\]
Let $U=\bigcup_{y\in Y}C(y)\subseteq A_{\ell+1}$.
Suppose to the contrary that
\begin{align}\label{supp U<Y}
|U|<|Y|=|\mathcal{F}_k|.
\end{align}
Note that $|Y|\le|W_{\ell+1}|\le n\le \xi N$, we have $e_{\Gamma}(\mathcal{F}_k,U)< 6p^k|U||\mathcal{F}_k|\le 6\xi Np^k|\mathcal{F}_k|$
since $G_r(N)\in\mathscr{C}_{N,p}^k(\xi)$ from property (ii) of Lemma \ref{main-le}.
However,
the fact that $\varphi_\ell(K(y))$ and $C(y)$ form a complete bipartite graph in $G_R$, which together with (\ref{Cy}) yield
\[
e_{\Gamma}(\mathcal{F}_k,U)\ge \frac{1}{4^{k+1}}p^k\frac{N}{T_0}|\mathcal{F}_k|\stackrel{\mathrm{(\ref{xi-ga-B})}}{\ge} 6\xi Np^k|\mathcal{F}_k|,
\]
which is a contradiction. This completes the proof of the claim. \hfill$\Box$

\smallskip
From Claim 2, we can extend $\varphi_{\ell}$ as follows:
\begin{align*}\label{de-varphi}
\varphi_{\ell+1}(w) = \left\{ \begin{array}{ll}
\varphi_{\ell}(w), & \textrm{if $w\in\bigcup_{j=1}^\ell W_j$,}\\
\psi(w), & \textrm{if $w\in  W_{\ell+1}$.}
\end{array} \right.
\end{align*}
From the inductive assumption of ($\mathcal{S}_{\ell}$) that $\varphi_{\ell}$ is an embedding of $H[\cup_{j=1}^\ell W_j]$ into $G_R[\cup_{j=1}^\ell A_j]$
 and $\psi$ is injective, we have $\varphi_{\ell+1}$ is indeed a partial embedding of $H[\bigcup_{j=1}^{\ell+1}W_j]$ into $G_R[\cup_{j=1}^{\ell+1} A_j]$.

Recall that $|N_H(z)\cap W_{\ell+1}|\le1$ for every $z\in \bigcup_{j=\ell+2}^{\bar{\Delta}}W_j$, so we define
\begin{align*}
C_{\ell+1}(z)=\left\{\begin{array}{ll}
C_{\ell}(z), &\textrm{if $N_H(z)\cap W_{\ell+1}=\emptyset$,}\\
C_{\ell}(z)\cap N_{G_R}(\varphi_{\ell+1}(y)), &\textrm{if $N_H(z)\cap W_{\ell+1}=\{y\}$.}
\end{array} \right.
\end{align*}
To complete the inductive step, it suffices to verify ($a$), ($b$) and ($c$) of
($\mathcal{S}_{\ell+1}$) for every $z\in \bigcup_{j=\ell+2}^{\bar{\Delta}}W_j$.

Fix  $z\in \bigcup_{j=\ell+2}^{\bar{\Delta}}W_j$.
If $N_H(z)\cap W_{\ell+1}=\emptyset$, then properties ($a$) and ($b$) of ($\mathcal{S}_{\ell+1}$) follows since $C_{\ell+1}(z)=C_{\ell}(z)$ and
 $\mbox{ldeg}_g^{\ell+1}(z)=\mbox{ldeg}_g^{\ell}(z)$. If $N_H(z)\cap W_{\ell+1}=\{y\}$, then property ($a$) of ($\mathcal{S}_{\ell+1}$) follows since $C_{\ell+1}(z)=C_{\ell}(z)\cap N_{G_R}(\varphi_{\ell+1}(y))$, and property ($b$) of ($\mathcal{S}_{\ell+1}$) follows from ($b'$) of Claim 1.

For property ($c$) of ($\mathcal{S}_{\ell+1}$), let $\{z,z'\}$ be an edge of $H$ with $z,z'\in \bigcup_{j=\ell+2}^{\bar{\Delta}}W_j$.
There are four cases according to the size of $N_H(z)\cap W_{\ell+1}$ and $N_H(z')\cap W_{\ell+1}$.
If $N_H(z)\cap W_{\ell+1}=\emptyset$ and $N_H(z')\cap W_{\ell+1}=\emptyset$,
then property ($c$) of ($\mathcal{S}_{\ell+1}$) follows directly from property ($c$) of ($\mathcal{S}_{\ell}$) and $\epsilon_{\ell+1}\ge \epsilon_{\ell}$.
If $N_H(z)\cap W_{\ell+1}=\{y\}$ and $N_H(z')\cap W_{\ell+1}=\emptyset$,
then property ($c$) of ($\mathcal{S}_{\ell+1}$) follows from ($c'$) of Claim 1; 
the case that $N_H(z)\cap W_{\ell+1}=\emptyset$ and $N_H(z')\cap W_{\ell+1}=\{y'\}$ is similar.
If $N_H(z)\cap W_{\ell+1}=\{y\}$ and $N_H(z')\cap W_{\ell+1}=\{y'\}$,
then $y=y'$, as otherwise $e_H(N_H(y), N_H(y'))\ge1$ will contradict (\ref{NwNw'=0}). Consequently, property ($c$) of ($\mathcal{S}_{\ell+1}$) follows from ($c'$) of Claim 1.

In conclusion, this completes the induction step and hence the proof of Lemma \ref{main-le}.\hfill$\Box$

\end{document}